\documentclass[fleqn,a4paper,12pt]{article}

\usepackage{amsmath}
\usepackage{amssymb}
\usepackage{bm}
\usepackage{amsthm}
\usepackage{graphicx}
\usepackage[round,sectionbib]{natbib}
\usepackage{setspace} 
\usepackage{a4wide}
\usepackage{url}
\usepackage{caption}
\usepackage{enumitem}
\usepackage{titling}
\usepackage{hyperref} 
\input{ee.sty}

\newtheorem{prp}{Proposition}

\newtheorem{con}{Condition}
\newtheorem{rmk}{Remark}

\thanksmarkseries{alph}

\renewcommand*{\eqref}[1]{%
	\hyperref[{#1}]{\textup{\tagform@{\ref*{#1}}}}%
}

\DeclareMathOperator*{\argmin}{arg\,min}

\allowdisplaybreaks
\begin{document}
	\begin{titlepage}
		\title{Exact Testing of Many Moment Inequalities\\ Against Multiple Violations}
		\vspace{1cm}
		\author{Nick Koning\thanks{
					Corresponding author. Nick Koning, Faculty of Economics and Business, University of Groningen, PO Box 800, 9700 AV Groningen, The Netherlands. Email: n.w.koning@rug.nl.
				}
				\hspace{0.035cm} and
				Paul Bekker\thanks{
					Paul Bekker, Faculty of Economics and Business, University of Groningen, Postbus 800, 9700 AV Groningen, The Netherlands. Email: p.a.bekker@rug.nl.
				}  \\
			{\small Faculty of Economics and Business} \\ {\small University of Groningen}
			\\}
		\date{\small May 27, 2020}
		
		\maketitle
		
		\thispagestyle{empty}
		\mbox{}\vspace{1cm}
		\begin{abstract}
			\noindent This paper considers the problem of testing many moment inequalities, where the number of moment inequalities ($p$) is possibly larger than the sample size ($n$). \citet{chernozhukov2019inference} proposed asymptotic tests for this problem using the maximum $t$ statistic. We observe that such tests can have low power if multiple inequalities are violated. As an alternative, we propose novel randomization tests based on a maximum non-negatively weighted combination of $t$ statistics. We provide a condition guaranteeing size control in large samples. Simulations show that the tests control  size in small samples ($n = 30$, $p = 1000$), and often has substantially higher power against alternatives with multiple violations than tests based on the maximum $t$ statistic.
			
		\end{abstract}
		
		\vspace{4cm}
		\noindent
		{\small
		}
		
	\end{titlepage}
	\newpage
	\doublespacing
	\section{Introduction}
		As discussed by \citet{chernozhukov2019inference}, henceforth CCK, the moments inequalities framework has developed into a powerful tool for inference on causal and structural parameters in partially identified models. In such models, the parameters of interest may be restricted to a subset of the parameter space defined by a collection of moment inequalities. The simultaneous testing of these moment inequalities provides inference about the true underlying parameter values. CCK provide an excellent review of the literature with detailed motivating examples. 
		
		They point out that many economic models give rise to problems where the number of moment inequalities $p$ may be much larger than the number of observations $n$. While there exists a large literature on testing moment inequalities, traditional methods are not well equipped for dealing with \textit{many} moment inequalities.\footnote{For work on unconditional moment inequalities see, e.g., \citet{canay2010inference, andrews2012inference, andrews2009validity, chernozhukov2007estimation, rosen2008confidence, romano2008inference}. \citet{andrews2013inference} notice that conditional moment inequalities can be viewed as an infinite number of unconditional moment inequalities. Contributions to conditional moment inequalities are found in \citet{chernozhukov2013intersection, lee2013testing, lee2018testing, armstrong2014weighted, armstrong2015asymptotically, armstrong2016multiscale}.} In order to test the moment inequalities against the alternative where at least one of them is violated, CCK use the maximum of $p$ $t$ values as a test statistic. They find critical values by using asymptotic theory and bootstrap methods. Additionally, a first-stage inequality selection step is included to improve the power of their tests. \citet{allen2018testing} suggests a refinement of the selection step. \citet{bugni2016inference} consider a generalization of the same problem and use Lasso for the first-stage selection. 
		
		Our contribution to the framework of testing many moment inequalities is threefold. First, we propose two novel test statistics. Notice that the maximum of $p$ $t$ values is invariant to the size of the second largest $t$ value. If the alternative hypothesis allows for no more than one violation, inference based on the maximum may be powerful, but it would discard power against alternatives where multiple moment inequalities are violated. 
		
		In order to retain this power, we propose an enhanced maximum statistic, where we maximize over an artificially expanded set of $t$ values. This expanded set of $t$ values is constructed by adding an extra set of $t$ values to the set of $p$ $t$ values. Such extra $t$ values are created by taking a non-negatively weighted combination of sample means divided by the square-root of its combined sample variance. This results in a test statistic that is also large if such a weighted combination is large, and not just if any of the original $t$ values is large. Interestingly, this expanded set of $t$ values is still a set of $t$ values (though correlated by construction). Therefore, the results of CCK can also be applied to this new test statistic, if the set of extra $t$ values is finite.
		
		In addition, we consider adding all non-negative combinations, so that the set of extra $t$ values is infinite. This results in a test statistic that can be viewed as a one-sided version of Hotelling's $T^2$ statistic. An interesting feature of our statistic is that, unlike the Hotelling's $T^2$ statistic, it can also be used in high dimensional settings. This is because its existence condition is weaker than non-singularity of the sample covariance matrix \citep{koning2019directing}. This condition is related to the positive eigenvalue condition \citep{meinshausen2013sign, slawski2013non} or copositivity of the sample covariance matrix. Like the positive eigenvalue condition, it may be satisfied if many elements of the sample covariance matrix are positive. However, it is also satisfied if many of the sample means are negative. In a simulation study, we find substantial increases in power compared to the statistic proposed by CCK, against alternatives where multiple moment inequalities are violated. However, a limitation of the test seems to be that it loses power against alternatives where the number of not satisfied moment inequalities is large compared to $n$ and the underlying covariance matrix exhibits weak correlations, as this interferes with the existence condition.
		
		Secondly, we propose using randomization tests in the moment inequalities framework, by imposing a symmetry assumption on the errors. Randomization tests originate with \citet{fisher1935design} and have been widely used in the literature.\footnote{See, for example, \citet{lehmann2006testing}, \citet{maritz1995distribution}, \citet{romano1990behavior}, \citet{bekker2002exact} and \citet{bekker2008symmetry}} An advantage of randomization tests is that they are exact if all moment inequalities are satisfied. In addition, we provide a sufficient condition to guarantee our tests control size for sufficiently large samples. This result also holds for the statistic that is the maximum over an infinite set of $t$ values.
		
		In our simulation experiments we find that the randomization tests perform similarly or slightly better compared to the empirical bootstrap proposed by CCK in large samples, and better in small samples. While randomization tests work well in small samples, they do require the validity of a symmetry assumption. However, in the simulation experiments where the symmetry assumption does not hold, we find, to our surprise, that the randomization test has better control of size than the empirical bootstrap procedure proposed by CCK, even if the sample is large. 
		
		Finally, we describe a first stage selection step in order to improve power by eliminating non-binding moment inequalities. The approach is based on symmetry and it does not depend on a specific selection method. Simulation results show that the tests with selection control size, and have increased power in the presence of non-binding moment inequalities compared to tests without selection. 
		
	\section{The model and two test statistics}\label{sec:teststat}
		Let observations in the $n\times p$ matrix $\mX$ satisfy
	\begin{align}\label{model}
	\mX = \viota_n\vmu' + \mE,
	\end{align}  where $\viota_n$ is an $n$-vector of ones, $\vmu$ is a $p$-vector of parameters and $\mE$ is a random matrix that satisfies a symmetry assumption. Specifically, we assume that any row of $\mE$ can be multiplied by $-1$ without affecting the distribution of $\mE$. This assumption is satisfied if the rows of $\mE$ are independently symmetrically distributed about zero. Moments need not exist and the rows of $\mE$ need not be identically distributed. We describe exact inference based on the finite sample.

	Following CCK, we are interested in testing the hypothesis H$_0$:\ $\vmu \leq \vzeros$ against the alternative H$_1$:\ $\vmu \not\leq \vzeros$.
	That is, we want to test whether \text{all} elements of $\vmu$ are non-positive, against the alternative that at least one is positive. We will refer to the null hypothesis as the \textit{moment inequalities}. The $j$th moment inequality is said to be \textit{violated} if $\mu_j > 0$, it is \textit{binding} if $\mu_j = 0$, and it is \textit{strictly satisfied} if $\mu_j < 0$.

	To define the test statistics, let $\mI_n = (\ve_1, \dots, \ve_n)$ be the $n \times n$ identity matrix, and define $\mP_{\viota_n} = \viota_n\viota_n'/\viota_n'\viota_n$. Define $\widehat{\vmu} = \mX'\viota_n/n$ with elements $\widehat{\mu}_j$, and $\widehat{\mSigma} = \mX'(\mI_n - \mP_{\viota_n})\mX/n$ with diagonal elements $\widehat{\sigma}_j^2$. We assume $\widehat{\vmu}\nleq\vzeros$, otherwise we would not reject H$_0$, and consider test statistics of the form
	\begin{align*}
	T_\calU=n^{1/2}\max_{\vlambda \in \calU} \frac{\widehat{\vmu}'\vlambda}{\sqrt{\vlambda'\widehat\mSigma\vlambda}},
	\end{align*}
	where $\calU\subseteq\mathbb{R}^p$ and Condition \ref{con:1} holds \citep{koning2019directing}.
	\begin{con}\label{con:1}
		$\vlambda'\widehat{\mSigma}\vlambda > 0$, for all $\vlambda \in\calU$, $\vlambda \neq \vzeros$,  $\widehat{\vmu}'\vlambda > 0$.
	\end{con} 
	\hspace{-.65cm}The test statistic has been used before for particular choices of $\calU$. Here we propose new choices as well.

	In order to test H$_0$  against H$_1$, CCK use  the maximum $t$ value, $t_j =n^{1/2} \widehat{\mu}_j/\widehat{\sigma}_j$, for $j =1, \dots, p$. We will denote this statistic by
	$
	t_{\text{max}}=t_{\text{max}}(\mX)
	= \max_{1 \leq j \leq p} t_j = T_{\calC^p}
	$, where $\calC^p=\{\ve_1, \dots, \ve_p\}$.
	By ignoring the sizes of all but the largest $t$ value, the $t_{\text{max}}$ statistic may  perform well against alternatives with one or few violations, but lack power when testing against alternatives with multiple violations. To shift power towards multiple violations, one might expand $\calU$ by adding vectors from the non-negative orthant, $\va_i\in\mathbb{R}_+^p$, $i=1,\ldots,m$. Let $\calU=\calC^p\cup\{\va_1,\ldots,\va_m\}$, then  $T_\calU$ may outperform $t_{\text{max}}$ in terms of power under alternatives with multiple violations.

	\begin{rmk}\label{rmk:1}
		Interestingly, if $\calU$ is finite, the $T_\calU$ statistic fits entirely within the framework of CCK. This implies that their results regarding the bootstrap and self-normalized methods can be immediately applied if the suitable regularity conditions hold. 
	\end{rmk}
	
	\begin{rmk}\label{rmk:2}
		If the off-diagonals of $\widehat{\mSigma}$ are positive, expanding the set $\ \calU$ may be less effective. In the extreme case that $\widehat{\mSigma} = \viota_p\viota_p'$, the $T_\calU$ and  $t_{\emph{max}}$ statistics coincide if $\calU\subset\mathbb{R}_+^p$ is finite and $\calC^p\subset\calU$. That is to say, $T_\calU \geq t_{\emph{max}}$ by construction, and			
		\begin{align*}
		\hspace{-.5cm}
		n^{-1/2}T_\calU
		= \max_{\lambda \in {\calU}} \frac{\widehat{\vmu}'\vlambda}{\sqrt{\vlambda'\widehat{\mSigma}\vlambda}}
		= \max_{\lambda \in{\calU}} \frac{\widehat{\vmu}'\vlambda}{\viota_p'\vlambda}
		= \max_{\substack{\lambda \in {\calU} \\ \viota_p'\vlambda = 1}}\widehat{\vmu}'\vlambda 
		\leq \max_{\substack{\lambda \geq \vzeros \\ \viota_p'\vlambda = 1}} \widehat{\vmu}'\vlambda
		= \max_{\lambda \in \calC^p}\widehat{\vmu}'\vlambda
		= n^{-1/2}t_{\emph{max}}.
		\end{align*}
	\end{rmk}

	Alternatively, we consider infinite sets $\calU$. In particular, we use $\calU=\mathbb{R}_+^p$ and define
	$
	T_+=T_{\mathbb{R}_+^p}.
	$\footnote{To handle cases with multiple violations, CCK mention the sum statistic $ \sum_{i = 1}^{p} (\sqrt{n}\max\{\widehat{\mu}_j/\widehat{\sigma}_j, 0\})^2$ and the quasi likelihood ratio test statistic $\min_{\vmu \leq \vzeros} n (\widehat{\vmu} - \vmu)'\widehat{\mSigma}^{-1}(\widehat{\vmu} - \vmu)$ as posibilities.	Only if $\widehat\mSigma$ is diagonal, these statistics fit within the $T_\calU$ framework. In that case the statistics are equal to $T_+$. However,	the sum statistic  ignores the covariances, and  the quasi likelihood-ratio statistic requires $\widehat{\mSigma}$ to be invertible.}
	The $T_+$ statistic can be viewed as a one-sided version of the square-root of Hotelling's $T^2$-statistic, since $T_{\mathbb{R}^p}^2=T^2 = n\widehat{\vmu}'\widehat{\mSigma}^{-1}\widehat{\vmu}$, which requires $\widehat\mSigma$ to be invertible. For the $T_+$ statistic $\widehat\mSigma$ need not be invertible. Instead,  a weaker condition suffices, which is a special case of  Condition \ref{con:1}.
	\begin{con}\label{con:2}
		$\vlambda'\widehat{\mSigma}\vlambda > 0$, for all $\vlambda \geq \vzeros$, $\vlambda \neq \vzeros$, $\widehat{\vmu}'\vlambda > 0$.
	\end{con}		
	\hspace{-.65cm}This condition is weaker than (strict) copositivity of $\widehat{\mSigma}$, which requires that $\vlambda'\widehat{\mSigma}\vlambda > 0$ for all $\vlambda \geq \vzeros$, $\vlambda \neq \vzeros$ (see e.g. \citealp{hiriart2010variational}), or the equivalent positive eigenvalue condition formulated by \citet{meinshausen2013sign}.\footnote{The positive eigenvalue condition $\min_{{\vlambda \geq \vzeros,\ \vlambda \neq \vzeros}} {\vlambda'\widehat{\mSigma}\vlambda}/{\|\vlambda\|_1^2} > 0$ is equivalent to  $\vlambda'\widehat{\mSigma}\vlambda > 0$, for all $\vlambda \geq \vzeros$, $\vlambda \neq \vzeros$.}
	
	An example where copositivity holds is if all elements of $\widehat{\mSigma}$ are strictly positive (i.e. $\widehat{\mSigma}$ is \textit{positive}). \citet{meinshausen2013sign} also provides an example where $\widehat{\mSigma}$ is allowed to contain some negative entries. In particular, if $\calA^c$ contains the indices of the the largest positive principal submatrix of $\widehat{\mSigma}$ and $\calA$ is its complement, then the condition only has to hold on the principal submatrix of $\widehat{\mSigma}$ with indices in $\calA$. If the number of elements in $\calA$ is small compared to $n$, this need not be very restrictive.
	
	The difference between Condition \ref{con:2} and copositivity lies in the role played by $\widehat{\vmu}$. Condition \ref{con:2} is equivalent to copositivity if $\widehat{\vmu} > \vzeros$. At the other extreme, if $\widehat{\vmu} \leq \vzeros$ the set $\{\vlambda\ |\ \vlambda \geq\vzeros,\ \lambda\neq\vzeros,\ \vmu'\vlambda\geq 0\}$ is empty, so that Condition \ref{con:2} holds, whether or not  copositivity holds.  If the elements of $\widehat{\vmu}$ have comparable positive or negative sizes, then Condition \ref{con:2} tends to become less restrictive the more negative elements there are. So, intuitively, even if $\widehat{\mSigma}$ is not copositive, Condition \ref{con:2} may still hold, in particular if many moment inequalities are strictly satisfied. This reasoning is confirmed in the simulation experiments, where the test based on $T_+$ performs well even when $n = 30$ and $p = 1000$, and only 10 moment inequalities are violated.
	
	To compute $T_+$ we use a recent result by \citet{koning2019directing}, who shows 
	\begin{align*}
	T_+&=n^{1/2}\max_{\vlambda \geq \vzeros} \frac{\widehat{\vmu}'\vlambda}{\sqrt{\vlambda'\widehat\mSigma\vlambda}}=n^{1/2}\frac{\widehat{\vmu}'\widehat\vlambda}{\sqrt{\widehat\vlambda'\widehat\mSigma\widehat\vlambda}},\quad\mbox{where}\quad 			\widehat\vlambda = \argmin_{\vlambda \geq \vzeros} \|\viota_n - \mX\vlambda\|_2^2.
	\end{align*}
	This problem can be solved by the non-negative least squares algorithm of \citet{lawson1995solving}.\footnote{This algorithm is implemented in the R package `nnls' and the MATLAB function `lsqnonneg'.} An interesting additional feature is that the resulting maximizing weight vector $\widehat{\vlambda}$ may contain zeros, resulting in a selection of moment inequalities that are `suspected' to violate the null.
			
	\section{Symmetry based inference}\label{sec:refl}
		The randomization tests that we consider are based on the symmetry assumption. For completeness, we include this section explain how the  assumption leads to exact randomization tests, and show how they can be applied in practice. For other descriptions of randomization tests, see e.g. \citet{lehmann2006testing}, \citet{maritz1995distribution} and \citet{bekker2008symmetry}. 

			Let $\calR = \{\mR_1, \mR_2, \dots, \mR_N\}$ be the set of $N = 2^n$ diagonal $n \times n$ matrices with diagonal elements in $\{-1, 1\}$, where $\mR_1$ denotes the identity matrix. This set constitutes a finite reflection group under matrix multiplication. The group $\calR$ determines a partitioning $\calE$ of $\mathbb{R}^{n \times p}$ into equivalence classes (``orbits") denoted by $\calR_{\mE^*} = \{\mE^*, \mR_2\mE^*, \dots, \mR_N\mE^*\}$, where $\mE^* \in \mathbb{R}^{n \times p}$ acts as a representative of the class. For simplicity, we only consider orbits $\calR_{\mE^*} \in \calE$, for which the cardinality satisfies $|\calR_{\mE^*}| = |\calR| = N$. This is equivalent to assuming that we only consider the subset of orbits $\calE^0 \subset \calE$ where the elements of $\calR_{\mE^*} \in \calE^0$ have no rows that equal zero.
			
			The basic assumption that permits the construction of randomization tests is a distributional invariance assumption on the error term under group transformations. In particular, we assume that $\mR\mE$ and $\mE$ have the same distribution, for all $\mR \in \calR$. Equivalently, we assume the conditional distribution of $\mE$, given an orbit  $\mE\in\calR_{\mE^*}$, is uniform for all orbits $\calR_{\mE^*} \in \calE^0$.

			Let $g: \calR_{\mE^*} \to \mathbb{R}\cup\{\infty\}$ be a mapping. Notice that if $g$ is injective on all orbits  $\calR_{\mE^*} \in \calE^0$, then the symmetry assumption implies that the conditional distribution of $g(\mE)$, given an orbit $\calR_{\mE^*} \in \calE^0$, is uniform just as well. Furthermore, the function $p(\mE) = |\{\mR \in \calR\ |\ g(\mR\mE) \geq g(\mE)\}| / N$ has a conditional distribution, given an orbit, that is uniform over $\{1/N, 2/N, \dots, 1\}$. As this holds for all $\calR_{\mE^*} \in \calE^0$, the function $p(\mE)$ has an unconditional uniform distribution  over $\{1/N, 2/N, \dots, 1\}$. This leads to the following results.

			\begin{prp}\label{prp:1}
				$\mbox{\emph{P}}(p(\mE) \leq \alpha) \leq \alpha$, for  $\alpha \in [0,1]$. If $g$ is injective on all $\calR_{\mE^*} \in \calE^0$, then $\mbox{\emph{P}}(p(\mE) \leq \alpha) = \alpha$, for $\alpha \in \{1/N, 2/N, \dots, 1\}$.
			\end{prp}					
			Following \citet{bekker2008symmetry}, we will refer to $g$ as an \textit{inferential function}. As the cardinality of $\calR$ grows exponentially in $n$, the computation of  $p(\mE)$ is intractable even in small samples. Therefore, we provide the following result to allow for sampling from $\calR$.\footnote{\citet{bekker2008symmetry} also describe sampling with replacement.}
			\begin{prp}\label{cor:without}
				Let $\calR^M$ contain the identity matrix and $M - 1$ other elements drawn randomly without replacement from $\calR \setminus \{\mR_1\}$. Let $p_M(\mE) = |\{\mR \in \calR^M\ |\ g(\mR\mE) \geq g(\mE)\}| / M$, then 
				$
				\mbox{\emph{P}}(p_M(\mE) \leq \alpha) \leq \alpha,
				$
				for $\alpha\in[0,1]$ and, if $g$ is injective, then 
								$\mbox{\emph{P}}(p_M(\mE) \leq \alpha) = \alpha$
								 for $\alpha \in\{1/M, 2/M, \dots, 1\}$.
			\end{prp}
			
			\section{Randomization tests for  moment inequalities}\label{sec:randtest}
			To test the moment inequalities of model \eqref{model}, we use randomization tests based on the $T_{\calU}$ statistic with $\calU \subseteq \mathbb{R}_+^p$, such as $t_{\text{max}}$ and $T_+$. We will use the convention that $T_{\calU} = 0$ if $\widehat{\vmu} \leq \vzeros$ and $T_{\calU} = \infty$ if Condition \ref{con:1} does not hold. Let the test statistic be generically denoted as $T$. The tests are exact if the moment inequalities are binding. In that case $\mX = \mE$ and $T$ may be seen as an inferential function $g(\mE) = T(\mE)$. So, Proposition \ref{cor:without}  can be used based on  $p_M(\mE) = p_M(\mX)=|\{\mR \in \calR^M\ |\  T(\mR(\mX))\geq T(\mX)| / M$.
			In particular, we use the following testing procedure. \\
			
			\noindent\textbf{Algorithm 1} (Symmetry Randomization Test).
			\begin{itemize}
				\item[1.] Create the set $\calR^M$ consisting of the $\mI_n$ and $M - 1$ matrices drawn from $\calR\setminus\{\mI_n\}$ without replacement.
				\item[2.] Generate the transformed dataset $\mR\mX$ and compute $T(\mR\mX)$, for all $\mR \in \calR_M$.
				\item[3.] Compute $p_M(\mX)$ by counting the proportion that exceeds $T(\mX)$.
				\item[4.] Reject the null if $p_M(\mX) \leq \alpha$.
			\end{itemize}
			If $\vmu=\vzeros$ and $T$ is injective and $\alpha\in\{1/N,2/N,\ldots,1\}$, this test is exact (note that point accumulation could occur at 0 or $\infty$). 

			In case some moment inequalities are assumed to be strictly satisfied, the test is no longer exact, but we observe in all our simulations that the size is controlled if the symmetry assumption holds. To prove this formally is another matter. If, under the null hypothesis $\vmu\leq\vzeros$, it can be verified that  $T(\mR\mX)\geq T(\viota_n\vmu'+\mR\mE)$, then the test is conservative.\footnote{Notice that for the numerator of  $T_+(\mR\mX)$ it holds that $\viota_n'\mR\mX\vlambda=\viota_n'\mR\viota_n\vmu'\vlambda+\viota_n'\mR\mE\vlambda\geq\viota_n'(\viota_n\vmu'+\mR\mE)\vlambda$.} In that case we find, that $T(\mX) = T(\viota_n\vmu' + \mE)$ is an inferential function, $g(\mE) = T(\viota\vmu' + \mE)$, so that
			\begin{align*}
				p^T(\mX)=|\{\mR \in \calR\ |\ T(\mR\mX) \geq T(\mX)\}|\ \geq\  |\{\mR \in \calR\ |\ {g}(\mR\mE) \geq {g}(\mE)\}|=p(\mE),
			\end{align*}
			which implies $\mbox{{P}}(p^T(\mX) \leq \alpha) \leq \mbox{{P}}(p(\mE) \leq \alpha) \leq \alpha$. Unfortunately, we have not been able to formally prove whether or not the condition holds for our test statistics. However, a condition that the test is conservative for sufficiently large samples is easy.
			
			\begin{con} \label{con:3}
				As $n \to \infty$ and $p$ is fixed, $\mE'\viota_n/n \overset{a.s.}{\to} \vzeros$ and $\mE'\mE/n \overset{a.s.}{\to} \mSigma$.
			\end{con}
			
			To verify $T(\mR\mX) \geq T(\viota_n\vmu' + \mR\mE)$, first notice that this inequality holds trivially if $\mR = \mI_n$ or $\vmu'\vlambda = 0$. So, we can restrict ourselves to the case that $\mR \neq \mI_n$ and $\vmu'\vlambda < 0$, as $\vmu'\vlambda \leq \vzeros$ for $\vlambda \geq \vzeros$. Furthermore, notice that 
			\begin{align*}
				T^*(\mX) 
					= \max_{\lambda \in {\calU}} \frac{\viota_n'\mX\vlambda}{\sqrt{\vlambda'\mX'\mX\vlambda}} 
					= \frac{T(\mX)}{\sqrt{1 + T(\mX)^2/n}},
			\end{align*}
			is a strict monotone transformation of $T$, so that we can instead look at $T^*$. Consider the function of $\vlambda \in \calU \subseteq \mathbb{R}_+^p$ that is maximized at $T^*(\mR\mX)$,
			
			\begin{align*}
			\frac{\viota_n'\mR\mX\vlambda}{\sqrt{\vlambda'\mX'\mX\vlambda}} 
			&= \frac{\viota_n'\mR\viota_n \vmu'\vlambda + \viota_n'\mR\mE\vlambda}{\sqrt{\vlambda'(\viota_n\vmu' + \mE)'(\viota_n\vmu' + \mE)\vlambda}}
			> \frac{n\vmu'\vlambda + \viota_n'\mR\mE\vlambda}{\sqrt{\vlambda'(\viota_n\vmu' + \mE)'(\viota_n\vmu' + \mE)\vlambda}} \\
			&= \frac{n\vmu'\vlambda + \viota_n'\mR\mE\vlambda}{\sqrt{n(\vmu'\vlambda)^2 + 2\vlambda'\mE'\viota_n\vmu'\vlambda + \vlambda'\mE'\mE\vlambda}},
			\end{align*}
			as $\viota_n'\mR\viota_n < n$. Consequently, due to Condition \ref{con:3}, as $\mE$ and $\mR\mE$ follow the same distribution,
			\begin{align*}
			n^{-1/2} \frac{\viota_n'\mR\mX\vlambda}{\sqrt{\vlambda'\mX'\mX\vlambda}}
			> \frac{\vmu'\vlambda}{(\vmu'\vlambda)^2 + \vlambda'\mSigma\vlambda} + o_{a.s.}(0).
			\end{align*}
			Furthermore, $T^*(\viota_n\vmu' + \mR\mE)$ maximizes
			\begin{align*}
			\frac{\viota_n'(\viota_n\vmu' + \mR\mE)\vlambda}{\sqrt{\vlambda'(\viota_n\vmu' + \mR\mE)'(\viota_n\vmu' + \mR\mE)\vlambda}}
			= \frac{n\vmu'\vlambda + \viota_n'\mR\mE\vlambda}{\sqrt{n(\vmu'\vlambda)^2 + 2\vlambda'\mE'\mR\viota_n\vmu'\vlambda + \vlambda'\mE'\mE\vlambda}},
			\end{align*}
			so that
			\begin{align*}
			n^{-1/2} \frac{\viota_n'(\viota_n\vmu' + \mR\mE)\vlambda}{\sqrt{\vlambda'(\viota_n\vmu' + \mR\mE)'(\viota_n\vmu' + \mR\mE)\vlambda}} 
			= \frac{\vmu'\vlambda}{(\vmu'\vlambda)^2 + \vlambda'\mSigma\vlambda} + o_{a.s.}(0).
			\end{align*}
			Let $\vlambda^*$ be the maximizer of $T^*(\viota_n\vmu' + \mR\mE)$, then
			\begin{align*}
			n^{1/2} T^*(\mR\mX) 
			\geq n^{-1/2} \frac{\viota_n'\mR\mX\vlambda^*}{\sqrt{\vlambda^{*\prime}\mX'\mX\vlambda^*}}
			> n^{-1/2} T^*(\viota_n\vmu' + \mR\mE) + o_{a.s.}(0).
			\end{align*}
			As a result, the condition $T^*(\mR\mX) \geq T^*(\viota_n\vmu' + \mR\mE)$ holds, and therefore $T(\mR\mX) \geq T(\viota_n\vmu' + \mR\mE)$, so the test is conservative if $n$ is sufficiently large.

	\section{Symmetry based inference with pre-selection} \label{sec:ineqsel}
		
		As the tests are found to be conservative in case some moment equalities are strictly satisfied, there may be room for power improvements. CCK  propose  an inequality selection step aiming at removal of such strict moment inequalities $j$. In particular, they remove moment inequalities $j$ for which $t_j < c$, where $c$ is some chosen cut-off value. They propose several techniques to select this cut-off value such that the asymptotic testing procedures remain applicable, up to a small correction of the significance level. A similar way to select the cut-off value using Lasso is proposed by \citet{bugni2016inference}.
		
		Our approach for inference with pre-selection is similar to the approach of the previous section where there was no pre-selection. We use symmetry based inference. Given any inequality selection rule we describe tests that are exact or conservative if the moment inequalities are binding. That is to say, if $\vmu=\vzeros$, then $T(\mX)=T(\mE)$ and  $p^T(\mX)=|\{\mR \in \calR\ |\ T(\mR\mX) \geq T(\mX)\}|\leq\alpha$.
		
		Let $\calJ(\mX)\subset\{1,\ldots, p\}$ be an index selection subset, and let $\mX_{\calJ(\mX)}$ be the submatrix of $\mX$ consisting of columns with indexes in $\calJ(\mX)$. As the selection depends on $\mX$, it is not guaranteed that
		$p^T(\mX_{\calJ(\mX)})=|\{\mR \in \calR\ |\ T(\mR\mX_{\calJ(\mX)}) \geq T(\mX_{\calJ(\mX)})\}|\leq\alpha$. However, if $\vmu=\vzeros$, then $T(\mX_{\calJ(\mX)})=T(\mE_{\calJ(\mE)})=g(\mE)$ is an inferential function and 
		\begin{align*}
		p_{\text{sel}}^T(\mX_{\calJ(\mX)})&=|\{\mR \in \calR\ |\ T(\mR\mX_{\calJ(\mR\mX)}) \geq T(\mX_{\calJ(\mX)})\}|\\
		&=|\{\mR \in \calR\ |\ {g}(\mR\mE) \geq {g}(\mE)\}|=p(\mE),
		\end{align*}
		which implies $\mbox{\emph{P}}(p_{\text{sel}}^T(\mX_{\calJ(\mX)})\leq \alpha)=\mbox{\emph{P}}(p(\mE) \leq \alpha) \leq \alpha$, for $\alpha\in[0,1]$ and, if $g$ is injective, then $\mbox{\emph{P}}(p_M(\mE) \leq \alpha) = \alpha$ for $\alpha \in\{1/M, 2/M, \dots, 1\}$, as in Proposition \ref{cor:without}. Therefore, we use the following testing procedure. \\
		\noindent\textbf{Algorithm 2} (Symmetry randomization test with inequality selection).
		\begin{itemize}
			\item[1.] Create the set $\calR^M$ consisting of the $\mI_n$ and $M - 1$ matrices drawn from $\calR\setminus\{\mI_n\}$ without replacement.
			\item[2.] Compute the selected inequalities  $\calJ(\mR\mX)$, for each $\mR \in \calR_M$.
			\item[3.] Compute $T(\mR\mX_{\calJ(\mR\mX)})$ for all $\mR \in \calR_M$.
			\item[4.] Compute  $p_{\text{sel}}^T(\mX)$, by counting the proportion that exceeds $T(\mX_{\calJ(\mX)})$.
			\item[5.] Reject the null if $p_{\text{sel}}^T(\mX) \leq \alpha$.
		\end{itemize}
		If $\vmu\leq\vzeros$ and $\vmu\neq\vzeros$, the test is conservative if $T(\mR\mX_{\calJ(\mR\mX)})\geq T\left((\viota\vmu'+\mR\mE)_{\calJ(\viota\vmu'+\mR\mE)}\right)$. The power of the test varies with the value of $\vmu$ and the selection method. We follow CCK and use $\calJ(\mX)=\left\{j\in\{1,\ldots,p\}\ |\ t_j>c\right\}$, where the selection constant is chosen using their empirical bootstrap procedure. The simulations show that the size of the tests is controlled when $\vmu \leq \vzeros$. 
		
	\section{Simulations}
		In this section, we present Monte Carlo simulation results. We use setups based on the simulation experiments presented by CCK and \citet{bugni2016inference}. The tests described in this paper are compared to the Empirical Bootstrap (EB) tests described in CCK.\footnote{We only compare to the EB tests described in CCK, as they find that the EB tests perform similarly to the Multiplier Bootstrap tests and better than the self-normalized tests.} All experiments were implemented in R and code for the tests will be made available at \url{https://github.com/nickwkoning/}.  \\
		
		\noindent\textbf{Data generation} \\
		The data is created as follows: we generate $n \times p$ matrix $\mX = \viota_n\vmu' + \mE\mA$, where $\viota_n$ is an $n$-vector of ones, $n\times p$ matrix $\mE$ has i.i.d. elements drawn from a distribution $F$, and $\mA$ is defined such that $\mSigma = \mA'\mA$, where $\mSigma$ has elements $\sigma_{ij} = \rho^{|i-j|}$.
		
		The parameters we use are $n \in \{30, 400\}$, $p \in \{200, 500, 1000\}$ and $\rho \in \{0, 0.5, 0.9\}$. For the vector $\vmu$, we consider four different designs. The signs of the elements of $\vmu$ are chosen to represent four general cases: all moment inequalities are binding (Design 1), most moment inequalities are strictly satisfied and some are binding (Design 2), all moment inequalities are violated (Design 3), and some moment inequalities are violated while some are strictly satisfied (Design 4). The values of $\vmu$ for each combination of $n$ and $p$ were selected to ensure that the rejection probabilities are bounded away from 1.
		
\begin{figure}
		\centering
		\includegraphics[width=6cm]{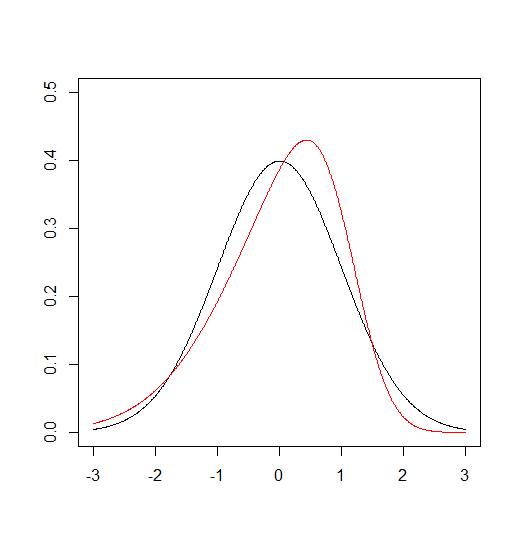}
		\includegraphics[width=6cm]{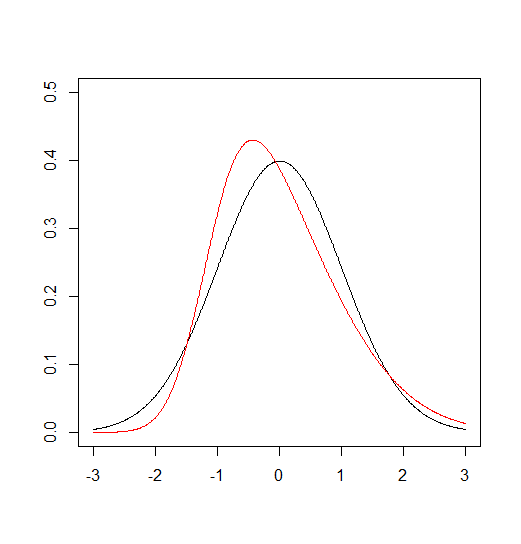}
		\caption{A density plot of the left-skewed (left) and right-skewed (right) normal distribution with mean 0, standard deviation 1 and skewness -0.667 and 0.667, respectively, overlaid on the density of the standard normal distribution.}
		\label{fig:skew}
	\end{figure}

		For the errors, we use a symmetric and two asymmetric distributions: one left-skewed and one right-skewed. As symmetric distribution, we use $t(4)/\sqrt{2}$, where $t(4)$ is the Student's $t$ distribution with 4 degrees of freedom and we divide by $\sqrt{2}$ so that the variance is 1. As asymmetric distributions we use the skew-normal distribution with mean 0, variance equal to 1 and two configurations for the skewness: $\gamma = -0.667$ and $\gamma = 0.667$. A density plot comparing these skewed distributions to the standard normal distribution is provided in Figure \ref{fig:skew}. For the sake of brevity, the asymmetric error distributions were only considered for Design 1.
		
		This setup is similar to the setup used by an earlier working paper of CCK.\footnote{This earlier version can be found at \url{https://arxiv.org/abs/1312.7614v4}.} The differences are that we do not consider equicorrelated data or uniformly distributed errors, but instead consider small samples ($n = 30$) and asymmetric error distributions. In addition, the positive values of $\vmu$ are substantially decreased to account for the higher power of tests based on the $T_+$ statistic. \\
		
		\noindent\textbf{Tests} \\
		We use symmetry based randomization (SR) tests based on the test statistics $t_{\text{max}}$ and $T_+$, with and without pre-selection. In addition, we also consider the $t_{\text{max}}^{\viota} = T_{\calC_p \cup \{\viota_p\}}$ statistic, where the maximum is taken over just a single extra vector compared to $t_{\text{max}}$. As a comparison, we include the EB tests described by CCK with and without pre-selection for both the $t_{\text{max}}$ and $t_{\text{max}}^{\viota}$ statistics.
		
		For each setting, 1000 realizations of $\mX$ were generated and the proportion of rejections was recorded. For the EB tests 1000 bootstrap samples were used and for the SR tests we used 1000 reflection samples. The significance level was fixed at $\alpha = 0.05$. The selection constant for the tests that include pre-selection is chosen using the EB procedure with $\beta = 0.001$ and 1000 resamples as described in CCK. For the small sample experiments with $n = 30$, we do not use pre-selection as the selection constant depends on the asymptotic properties of the EB test.
		
		In terms of computation time for the largest settings, a single test without pre-selection using the $T_+$ statistic a on single core of a standard 13-inch 2017 MacBook Pro takes approximately 80-270 seconds, depending on the design. For the $t_{\text{max}}$ and $t_{\text{max}}^{\viota}$ statistics this is approximately 5-7 seconds.
		
		\subsection{Results}
		The results of the simulation experiments can be found in Tables \ref{tab:1} to \ref{tab:5}, and are discussed separately for each design. \\
		
		\noindent\textbf{Design 1}: $\vmu = \vzeros$. \\
		In Design 1, all elements of $\vmu$ are equal to zero. Therefore, the rejection probability for the tests should be at most $\alpha = 0.05$. The results for Design 1 are displayed in Table \ref{tab:1} for the symmetric error distribution, and Table \ref{tab:5} for the asymmetric error distributions.
		
		In Table \ref{tab:1}, for the case that the number of observations is large ($n = 400$), we see that the EB tests for both $t_{\text{max}}$ and $t_{\text{max}}^{\viota}$ reject with probability close to 0.05, which is in line with Remark \ref{rmk:1}. For the symmetry randomization (SR) tests, Proposition \ref{cor:without} states that all rejection probabilities should be at most 0.05, which is confirmed in the simulations. In addition, if the test statistic is injective on all orbits, then the rejection probability is equal to $\alpha$. One exception found to this in the simulations is the configuration $p = 1000$ and $\rho = 0$, for the statistic $T_+$, where the proportion of rejections is 0. Further inspection shows that for this configuration without pre-selection, the test statistics $T_+(\mX)$ and critical values, defined as the $1 - \alpha$ quantile of the set $\{T_+(\mR\mX)\ |\ \mR \in \calR_M\}$, are infinite, as Condition \ref{con:2} is not satisfied. Therefore, the function $T_+$ is not injective on the orbit and the rejection probability can be strictly smaller than $\alpha$, according Proposition \ref{cor:without}. A similar problem of non-injectivity occurs for the case with pre-selection. In line with the relation between Condition \ref{con:2} and copositivity, the issue diminishes if $\rho$ is increased.
		
		If the number of observations is small ($n = 30$), then for the $t_{\text{max}}$ and $t_{\text{max}}^{\viota}$ statistics, the rejection rates for the SR tests remain approximately $\alpha$. In contrast, the EB tests over-reject. For the SR tests based on the $T_+$ statistic, the rejection rate is frequently zero due to the non-injectivity phenomenon described above.
		
		Table \ref{tab:5} shows the rejection rate for Design 1 when the error distributions are left-skewed ($\gamma = -0.667$) and right-skewed ($\gamma = 0.667$). As the symmetry assumption is violated, it is not surprising that the SR tests have a rejection rate different from $\alpha$. In particular, we find that the rejection rate for the left-skewed error distribution is larger than $\alpha$, while the rejection rate for the right-skewed error distribution is smaller than $\alpha$.
		
		Surprisingly, even though the sample size is large ($n = 400$), the EB tests over-reject even more than the SR tests under the left-skewed distribution. This suggests that the EB tests may also benefit from errors being symmetrically distributed in finite samples. Although these simulation results are by no means exhaustive, they suggest that one may sometimes be better off using an SR test than an EB test even if it is known that the true error distribution is not symmetric. \\
		
		\noindent\textbf{Design 2}: $\vmu \leq \vzeros$, $\vmu \neq \vzeros$. \\
		In Design 2, for $n = 400$, the first $0.1p$ elements of $\vmu$ are equal to 0 and the remaining elements are equal to $-0.8$. For $n = 30$, the first $10$ elements are 0 and the remaining elements are equal to $-5$. As the data for Design 2 is generated under $H_0$, the rejection rates should be smaller than $\alpha$ (up to sampling errors) if no pre-selection is used, and close to $\alpha$ if pre-selection is used. The results for Design 2 can be found in Table \ref{tab:2}.
		
		From Table \ref{tab:2}, it can be observed that the rejection rates of the tests without pre-selection is smaller than $\alpha$. For the SR tests without pre-selection and $n = 400$, this is in line with the conclusions in Section \ref{sec:randtest} as $n$ is large. For $n = 30$, the rejection rates are closer to, but still smaller than $\alpha$.  In addition, pre-selection lifts the rejection rates close to the nominal size $\alpha$. \\
		
		\noindent\textbf{Design 3}: $\vmu > \vzeros$. \\
		In Design 3, all elements of $\vmu$ are positive. In particular, for $n = 400$ they are equal to 0.01 and for $n = 30$ they are equal to 0.03. As the data is generated under the alternative, the rejection probability should be as large as possible. The results for Design 3 are found in Table \ref{tab:3}.
		
		In Table \ref{tab:3}, it can be seen that for $n = 400$ the power of the EB and SR tests is similar, and inequality selection has no noticeable effect. The $T_+$ statistic generally leads to the highest power. One exception is the case where $p = 1000$ and $\rho = 0$, which is the most challenging case with respect to Condition \ref{con:2}. The power difference between the $t_{\text{max}}^{\viota}$ and $t_{\text{max}}$ statistics is quite remarkable, as the $t_{\text{max}}^{\viota}$ statistic maximizes over just a single extra vector. However, in Design 4 we will see that the statistic is more sensitive to the alternative than $T_+$. Furthermore, the difference in the power between the statistics is largest under weak correlations, which is in line with Remark \ref{rmk:2}.
		
		For $n = 30$, the results for the EB tests should be ignored, as the tests do not control size. Although there are multiple violations, the $T_+$ statistic performs worst. This is expected as Condition \ref{con:2} is unlikely to be satisfied here. So the same phenomenon observed under Design 1 occurs, where both the test statistic and critical value are infinite. Therefore, it is recommended that the test statistic and critical value are inspected before the outcome of the test is interpreted. If both are infinite, then the test based on the $T_+$ statistic should not be used. In addition, the $t_{\text{max}}^{\viota}$ statistic generally performs as well as or better than the $t_{\text{max}}$ statistic, especially if $\rho = 0$. \\
		
		\noindent\textbf{Design 4}: $\vmu \not\leq \vzeros$, $\vmu \not > \vzeros$. \\
		In Design 4, for $n = 400$, the first $0.1p$ elements of $\vmu$ are equal to 0.02 and the remaining elements are equal to $-0.75$. For $n = 30$, the first 10 elements of $\vmu$ are equal to 0.3 and the remaining elements are equal to -5. The data is generated under the alternative, so the rejection rates should be as large as possible. The results for Design 4 are presented in Table \ref{tab:4}.
		
		The results in Table \ref{tab:4} for $n = 400$ show that the tests with pre-selection have substantially higher power than the tests without pre-selection. When considering the tests with pre-selection, tests based on the $T_+$ statistic have much more power than the tests based on the $t_{\text{max}}$ and $t_{\text{max}}^{\viota}$ statistics. The power for the EB and SR tests based on the $t_{\text{max}}$ statistic are similar if pre-selection is used, but the SR tests seem more powerful if pre-selection is not used. In addition, the performance of the $t_{\text{max}}^{\viota}$ and $t_{\text{max}}$ statistics is similar, unlike in Design 3. Further inspection shows that the pre-selection often fails to eliminate a few strictly satisfied moment inequalities, which causes the added weight vector $\viota_p$ to rarely be the maximizing vector.
		
		Even though the EB tests do not control size for $n = 30$, the SR tests seem to be more powerful. Comparing the different statistics for the SR tests, it can be observed that the $T_+$ statistic performs well compared to the $t_{\text{max}}$ statistic if $\rho = 0$, even if $p = 1000$ and $n = 30$. So, the issue regarding Condition \ref{con:1} is not observed here.  This is in line with the intuitions obtained in Section \ref{sec:teststat} which suggest that the condition may still be satisfied if many of the moment inequalities are strictly satisfied. The $t_{\text{max}}^{\viota}$ statistic performs best overall for $n = 30$.
	
	\singlespacing
	\bibliographystyle{abbrvnat}
	\bibliography{bibfile}
	\clearpage
	
	\begin{table}[ht]
		
		\centering
		\begin{tabular}{rrrcccc|cccc|cc}
			\hline
			\hline
			&&&\multicolumn{4}{c}{$t_{\text{max}}$} & \multicolumn{4}{c}{$t_{\text{max}}^{\viota}$}  & \multicolumn{2}{c}{$T_+$} \\
			\hline
			n & p & $\rho$ & EB & SR & EB & SR & \multicolumn{1}{c}{EB} & SR & \multicolumn{1}{c}{EB} & SR & \multicolumn{1}{c}{SR} & \multicolumn{1}{c}{SR}\\
			& &  &  & & sel & sel & &  & sel & sel && sel  \\ 
			\hline
			&		& 0   & .045 & .046 & .032 & .052 & .042 & .050& .041& .050 & .040 & .043\\ 
			& 200 	& 0.5 & .044 & .038 & .046 & .057 & .050 & .048& .042& .061 & .055 & .046\\ 
			&		& 0.9 & .047 & .049 & .052 & .040 & .052 & .051& .075& .060& .068 & .050\\ 
			\hline
			&		& 0   & .044 & .054 & .043 & .032 & .057 & .042 & .041&.051& .054 & .061 \\ 
			400 & 500 	& 0.5 & .043 & .054 & .034 & .069 & .051 & .049& .038& .049& .054 & .045\\ 
			&		& 0.9 & .059 & .051 & .053 & .046 & .060& .050 & .057& .059& .048 & .049\\ 
			\hline
			&		& 0   & .043 & .042 & .051 & .063 & .036 & .060& .036& .053& .000 & .000 \\ 
			&1000	& 0.5 & .044 & .063 & .051 & .052 & .052 & .057& .038& .052& .049 & .054\\ 
			&		& 0.9 & .058 & .057 & .042 & .052 & .046 & .040& .047& .041& .039 & .051\\ 
			\hline
			\hline
			&   	& 0  	& .094 & .055 &-&-& .111& .061& -& -&   .000 & -\\
			&200 	& 0.5 	& .122 & .056 &-&-& .117& .048& -& - & .000 & -\\
			&   	& 0.9 	& .128 & .050 &-&-& .153 & .045 &- & - & .048 & -\\
			\hline
			&		& 0 	& .114 & .040 &-&-& .124 & .045 &-& - & .000 & -\\
			30	&500	& 0.5 	& .142 & .044 &-&-& .141 & .049& -& - & .000 & -\\
			&		& 0.9 	& .167 & .052 &-&-& .173 & .058&- & - & .000 & -\\
			\hline 
			&		& 0		& .137 & .049 &-&-& .155& .045& -& - & .000 & -\\
			&1000	& 0.5 	& .170 & .063 &-&-& .174 & .051& - & - & .000 & -\\
			&		& 0.9 	& .210 & .049 &-&-& .204 & .060 &- & - & .000 & -\\
			\hline
			\hline& & 
		\end{tabular}
		\caption{Monte Carlo rejection probabilities with 1000 repetitions for Design 1: $\vmu = \vzeros$, with symmetrically distributed errors ($t(4)/\sqrt{2}$). The columns represent the Empirical Bootstrap (EB) and Symmetry Randomization (SR) tests, based on the $t_{\text{max}}$, $t_{\text{max}}^{\viota}$ and $T_+$ statistics, both with pre-selection (sel) and without pre-selection.}
		\label{tab:1}
	\end{table}
	\begin{table}
		\centering
		\begin{tabular}{rrrcccc|cccc|cc}
			\hline
			\hline
			&&&\multicolumn{4}{c}{$t_{\text{max}}$} & \multicolumn{4}{c}{$t_{\text{max}}^{\viota}$} & \multicolumn{2}{c}{$T_+$}\\
			\hline
			n & p & $\rho$ & EB & SR & EB & SR & EB & SR & EB & SR & SR & SR \\
			& &  &  & & sel & sel && & sel & sel & & sel \\ 
			\hline
			&		& 0   & .003 & .009 & .038 & .037 & .002&.009& .033&.053& .000 & .049\\ 
			&200 	& .5 & .005 & .006 & .031 & .045 & .006& .004& .045&.048& .001 & .055\\ 
			&		& .9 & .006 & .009 & .040 & .054 & .004&.011& .041 &.038& .000 & .060\\ 
			\hline
			&		& 0   & .000 & .011 & .035 & .057 & .003& .003& .043&.032& .002 & .042\\ 
			400	& 500	& .5 & .008 & .005 & .059 & .063 & .004&.008& .046 &.037& .001 & .047\\ 
			&		& .9 & .005 & .005 & .039 & .049 & .005& .009& .052 &.050& .000 & .048\\ 
			\hline
			&		& 0   & .002 & .009 & .025 & .050 & .004 &.003& .032 &.038& .000 & .057 \\ 
			& 1000 	& .5 & .004 & .012 & .040 & .053 & .003 &.007& .031 &.050& .000 & .055\\ 
			&		& .9 & .007 & .012 & .047 & .049 & .003 &.010& .041&.030& .001 & .050\\ 
			\hline
			\hline
			&		& 0   & .004 & .034  &-&-& .006&.044 &-& - &  .009 & -\\
			& 200 	& .5 & .006 & .020  &-&-& .013 &.046&-&- &  .012 & -\\
			&		& .9 & .010 & .023  &-&-& .014& .033&-&- &  .014 & -\\
			\hline
			&		& 0   & .002 & .021  &-&-&.011&.042&-&- &  .008 & -\\
			30 	& 500	& .5 & .005 & .028  &-&-& .019&.040&-& -&  .010 & -\\
			&		& .9 & .002 & .037  &-&-& .017& .045&-& - &  .007 & -\\
			\hline 
			&		& 0	  & .002 & .035  &-&-& .014&.037&-& - &  .003 & -\\
			& 1000	& .5 & .003 & .022  &-&-& .024&.032&-& -&  .013 & -\\
			&		& .9 & .004 & .024  &-&-& .029&.042&-&- &  .006 & -\\
			\hline
			\hline
		\end{tabular}
		\caption{Monte Carlo rejection probabilities with 1000 repetitions for Design 2:  $\vmu \leq \vzeros$, with symmetrically distributed errors ($t(4)/\sqrt{2}$). For the cases that $n = 400$, $\mu_j = 0$ if $j \leq 0.1p$ and $\mu_j = -0.8$ if $j > 0.1p$, and for the cases where $n = 30$, $\mu_j = 0$ if $j \leq 10$ and $\mu_j = -5$ if $j > 10$. The columns represent the Empirical Bootstrap (EB) and Symmetry Randomization (SR) tests,  based on the $t_{\text{max}}$, $t_{\text{max}}^{\viota}$ and $T_+$ statistics, both with pre-selection (sel) and without pre-selection.}
		\label{tab:2}
	\end{table}
	\begin{table}
		\centering
		\begin{tabular}{rrrcccc|cccc|cc}
			\hline
			\hline
			&&&\multicolumn{4}{c}{$t_{\text{max}}$} & \multicolumn{4}{c}{$t_{\text{max}}^{\viota}$}& \multicolumn{2}{c}{$T_+$}\\
			\hline
			n & p & $\rho$ & EB & SR & EB & SR & EB & SR & EB & SR &\multicolumn{1}{c}{SR} & \multicolumn{1}{c}{SR} \\
			& &  &  & & sel & sel & & & sel & sel & & \multicolumn{1}{c}{sel} \\ 
			\hline
			&		& 0   & .093 & .099 & .094 & .099 & .307 &.333& .313& .340& .668 & .667 \\
			& 200 	& .5 & .099 & .087 & .094 & .084 & .135 & .125 & .141 & .117&.377 & .357 \\ 
			&		& .9 & .099 & .119 & .099 & .099 & .086 & .114& .078& .083& .129 & .144 \\
			\hline
			&		& 0   & .103 & .108 & .094 & .112 & .755 & .806& .773 & .799& .919 & .926 \\
			400	& 500 	& .5 & .106 & .107 & .100 & .097 & .206& .216& .207 & .237& .649 & .653 \\
			&		& .9 & .116 & .083 & .111 & .094 & .092 & .109& .092& .108& .215 & .230 \\
			\hline
			&		& 0   & .095 & .106 & .094 & .100 & .985 & .991 & .991 &.994& .000 & .000 \\
			& 1000 	& .5 & .099 & .116 & .094 & .095 & .447&.458& .461 & .453& .848 & .864 \\
			&		& .9 & .103 & .117 & .111 & .099 & .113& .119& .129& .111& .365 & .342 \\
			\hline
			\hline
			&		& 0		& .189 & .107&-&-& .297& .202 & - & - & .000 &-\\
			& 200& .5	& .190 & .101 &-&-& .231& .094& - & - & .000 &-\\
			&		& .9	& .167 & .086 &-&-& .168& .075 & - & - & .129 &-\\
			\hline
			&		& 0   & .225 & .092 &-&-& .610 & .411 & - & - & .000 &-\\
			30	& 500	& .5 & .253 & .100 &-&-& .298&.121& - & - &.000 &-\\
			&		& .9 & .275 & .099 &-&-& .264& .068& - & - & .000 &-\\
			\hline 
			&		& 0   & .242 & .111 &-&-& .896& .818& - & - & .000 &-\\
			& 1000	& .5 & .302& .091 &-&-& .429 & .209& - & - &.000 &-\\
			&		& .9 & .350 & .086 &-&-& .315&.115& - & - & .000&-\\
			\hline
			\hline
		\end{tabular}
		\caption{Monte Carlo rejection probabilities with 1000 repetitions for Design 3: $\vmu > \vzeros$, with symmetrically distributed errors ($t(4)/\sqrt{2}$). For the cases that $n = 400$, $\mu_j = 0.01$ and for the cases that $n = 30$, $\mu_j = 0.03$, for all $j$. The columns represent the Empirical Bootstrap (EB) and Symmetry Randomization (SR) tests,  based on the $t_{\text{max}}$, $t_{\text{max}}^{\viota}$ and $T_+$ statistics, both with pre-selection (sel) and without pre-selection.}
		\label{tab:3}
	\end{table}
	\begin{table}
		\centering
		\begin{tabular}{rrrcccc|cccc|cc}
			\hline
			\hline
			&&&\multicolumn{4}{c}{$t_{\text{max}}$} & \multicolumn{4}{c}{$t_{\text{max}}^{\viota}$} & \multicolumn{2}{c}{$T_+$} \\
			\hline
			n & p & $\rho$ & EB & SR & EB & SR &EB & SR &EB & SR & SR & SR \\
			& &  &  & & sel & sel & & & sel & sel && sel \\ 
			\hline
			&& 0 & .017 & .033 & .135 & .133 &.018 & .034& .109 & .124& .025 & .394\\
			&200 & .5 & .026 & .032 & .136 & .139 & .010 & .035& .147 & .129& .006 & .242\\
			&& .9 & .020 & .027 & .113 & .122 &.011 &.034& .123& .111& .002 & .130 \\
			\hline
			&& 0 & .018 & .031 & .099 & .103 & .024& .029& .107& .136& .143 & .675\\
			400&500 & .5 & .025 & .038 & .141 & .136 & .022 & .034& .138& .156& .015& .367\\
			&& .9 & .026 & .037 & .117 & .146 & .015& .022& .137&.126 & .002 & .159 \\
			\hline
			&& 0 & .022& .049 & .091 & .119 &.020 & .039&.082&.131& .397 & .898\\
			&1000 & .5 & .022& .041 & .153 & .151 & .014& .032& .126& .164& .087 & .550\\
			&& .9 & .030 & .037 & .135 & .131 &.021 & .033& .140 &.143& .003 & .214\\
			\hline
			\hline
			&		& 0	  & .509 & .828 &-&-& .734& .925 & - & - & .976 &-\\
			& 200 	& .5 & .417 & .709 &-&-& .645& .846 & - & - & .721 &-\\
			&		& .9 & .277 & .486 &-&-& .421 & .571 & - & - & .319 &-\\
			\hline
			&		& 0   & .342 & .813 &-&-& .887& .982 & - & - &.966 &-\\
			30	& 500	& .5 & .298 & .701 &-&-& .810 & .937 & - & - & .673 &-\\
			&		& .9 & .224 & .499 &-&-& .550& .693 & - & - & .318 &-\\
			\hline 
			&		& 0   & .278 & .802 &-&-& .952 & .996 & - & - & .951 &-\\
			& 1000	& .5 & .237 & .687 &-&-& .930 & .980 & - & - & .642 &-\\
			&		& .9 & .156 & .466 &-&-& .653 & .799& - & - & .296 &-\\
			\hline
			\hline& & 
		\end{tabular}
		\caption{Monte Carlo rejection probabilities with 1000 repetitions for Design 4: $\vmu \not \leq \vzeros$, $\vmu \not>\vzeros$, with symmetrically distributed errors ($t(4)/\sqrt{2}$). For the cases that $n = 400$, $\mu_j = 0.02$ if $j \leq 0.1p$ and $\mu_j = -0.75$ if $j > 0.1p$, and for the cases that $n = 30$, $\mu_j = 0.3$ if $j \leq 10$ and $\mu_j = -5$ if $j > 10$. The columns represent the Empirical Bootstrap (EB) and Symmetry Randomization (SR) tests,  based on the $t_{\text{max}}$, $t_{\text{max}}^{\viota}$ and $T_+$ statistics, both with pre-selection (sel) and without pre-selection. }
		\label{tab:4}
	\end{table}
	\clearpage	
	
	\begin{table}
		\centering
		\begin{tabular}{rrrcccc|cccc|cc}
			\hline
			\hline
			&&&\multicolumn{4}{c}{$t_{\text{max}}$} &\multicolumn{4}{c}{$t_{\text{max}}^{\viota}$} & \multicolumn{2}{c}{$T_+$}\\
			\hline
			$\gamma$ & p & $\rho$ & EB & SR & EB & SR & EB & SR & EB & SR & \multicolumn{1}{c}{SR} & \multicolumn{1}{c}{SR} \\
			& &  &  & & sel & sel && & sel & sel & & \multicolumn{1}{c}{sel} \\ 
			\hline
			&	 	& .0 & .116 	& .069	& .104 & .089 & .086 & .077& .082& .106 & .082 & .081 \\
			&200	 	& .5 & .107 	& .083 	& .098 & .063 & .081 & .061& .072& .083 & .084 & .073 \\
			&	 	& .9 & .069 	& .053 	& .075 & .050 & .056 & .048& .057& .059 & .070 & .059 \\
			\hline
			&	 	& .0 & .119 	& .087 	& .115 & .093 & .124 & .084 & .081 & .135 & .081 & .098 \\
			-.667	&500  	& .5 & .120 	& .086 	& .097 & .082 & .102 & .058 & .080 & .093 & .074 & .090 \\
			&	 	& .9 & .064 	& .058 	& .074 & .071 & .060 & .056 & .058 & .059 & .056 & .048 \\
			\hline
			&	 	& .0 & .136 	& .087 	& .146 & .093 & .153 &.091& .106 & .136& .000 & .000 \\
			&1000 	& .5 & .109 	& .075 	& .112 & .093 & .109 &.104& .089 & .121 & .068 & .085 \\
			&	 	& .9 & .065 	& .058 	& .087 & .056 & .059 &.065& .074& .067& .054 & .052 \\
			\hline
			\hline
			&		& .0 & .031 	& .025 	& .034 & .024 & .017 & .034& .032& .020& .043 & .027 \\
			&200 	& .5 & .038 	& .034 	& .034 & .031 & .033 & .040 & .042& .042& .035 & .036 \\
			&		& .9 & .042 	& .037 	& .050 & .043 & .045 & .048 & .034 & .042& .045 & .053 \\
			\hline
			&		& .0 & .020 	& .019 	& .023 & .033 & .018 & .024 &.017& .024& .021 & .019 \\
			.667	&500 	& .5 & .033 	& .024 	& .037 & .041 & .026 & .036 &.027& .028& .023 & .027 \\
			&		& .9 & .030 	& .037 	& .039 & .042 &.049 & .048 &.054& .050& .039 & .047 \\
			\hline
			&		& .0 & .033 	& .025 	& .021 & .030 &.019 & .028 &.016& .027 & .000 & .000 \\
			&1000	& .5 & .029 	& .036 	& .029 & .021 &.034 & .031&.025& .028 & .038 & .034 \\
			&		& .9 & .059 	& .042 	& .036 & .048 & .043 & .052&.051& .041& .034 & .046 \\
			\hline
			\hline
		\end{tabular}
		\caption{Monte Carlo rejection probabilities with 1000 repetitions for Design 1: $\vmu = \vzeros$, $n = 400$ observations, asymmetrically distributed errors (skew-normal with skewness parameter $\gamma$). The columns represent the Empirical Bootstrap (EB) and Symmetry Randomization (SR) tests,  based on the $t_{\text{max}}$, $t_{\text{max}}^{\viota}$ and $T_+$ statistics, both with pre-selection (sel) and without pre-selection.}
		\label{tab:5}
	\end{table}
	
\end{document}